\documentclass[12pt,reqno]{amsart}

\usepackage{amsfonts, amsthm, amsmath}
\allowdisplaybreaks[4]

\usepackage{rotating}

\usepackage{tikz}

\usepackage{graphics}

\usepackage{amssymb}

\usepackage{amscd}

\usepackage[latin2]{inputenc}

\usepackage{t1enc}

\usepackage[mathscr]{eucal}

\usepackage{indentfirst}

\usepackage{graphicx}

\usepackage{graphics}

\usepackage{pict2e}

\usepackage{mathrsfs}

\usepackage{enumerate}
\usepackage[pagebackref]{hyperref}
\hypersetup{backref, pagebackref, colorlinks=true}
\usepackage{cite}
\usepackage{color}
\usepackage{epic}
\usepackage{hyperref} 
\numberwithin{equation}{section}
\theoremstyle{plain}

\newtheorem{theorem}{Theorem}[section]

\newtheorem{lemma}[theorem]{Lemma}

\theoremstyle{definition}

\newtheorem{?}[theorem]{Problem}

\def\boxit#1{\leavevmode\hbox{\vrule\vtop{\vbox{\kern.33333pt\hrule
    \kern1pt\hbox{\kern1pt\vbox{#1}\kern1pt}}\kern1pt\hrule}\vrule}}

\usepackage{collectbox}



\newcommand{\f}[1]{\ifthenelse{\equal{#1}{1}}{(q;q)_\infty}{(q^{#1};q^{#1})_{\infty}}}

\begin{document}

\title[New congruences for broken $k$-diamond partitions]{New congruences for broken $k$-diamond partitions}

\author[D. Tang]{Dazhao Tang}

\address[Dazhao Tang]{College of Mathematics and Statistics, Chongqing University, Huxi Campus LD206, Chongqing 401331, P.R. China}
\email{dazhaotang@sina.com}

\date{\today}

\begin{abstract}
The notion of broken $k$-diamond partitions was introduced by Andrews and Paule. Let $\Delta_{k}(n)$ denote the number of broken $k$-diamond partitions of $n$ for a
fixed positive integer $k$. In this paper, we establish new infinite families of broken $k$-diamond partition congruences.
\end{abstract}

\subjclass[2010]{05A17, 11P83}

\keywords{Partitions; Broken $k$-diamonds; Congruences}

\maketitle



\section{Introduction}\label{sec1}

In 2007, Andrews and Paule \cite{AP2007} introduced a new class of directed graphs, called broken $k$-diamond partitions. They proved that the generating function of $\Delta_{k}(n)$, the number of broken $k$-diamond partitions of $n$, is given by
\begin{align*}
\sum_{n=0}^{\infty}\Delta_{k}(n)q^{n}=\dfrac{\f{2}\f{2k+1}}{\f{1}^3\f{2(2k+1)}}.
\end{align*}

Here and in the sequel, we assume $|q|<1$ and adopt the following customary notation on $q$-series and partitions:
\begin{align*}
(a;q)_{\infty} &:=\prod_{n=0}^{\infty}(1-aq^{n}).
\end{align*}

The following two congruences modulo 5 were subsequently proved by Chan \cite{Chan2008}, Radu \cite{Rad2015}, and Hirschhorn \cite{Hir2017}:
\begin{align*}
\Delta_{2}(25n+14) &\equiv0\pmod{5},\\
\Delta_{2}(25n+24) &\equiv0\pmod{5}.
\end{align*}

Moreover, a number of infinite families of congruences modulo 5 satisfied by $\Delta_{2}(n)$ have been proved. See, for example, Chan \cite{Chan2008}, Hirschhorn \cite{Hir2017}, Radu \cite{Rad2015} and Xia\cite{Xia2017}.

On the other hand, Jameson\cite{Jam2013} and Xia\cite{Xia2015} proved the following congruences modulo 7 enjoyed by $\Delta_{3}(n)$, which was conjectured by Paule and Radu \cite{PR2010}:
\begin{align*}
\Delta_{3}(343n+82) &\equiv0\pmod{7},\\
\Delta_{3}(343n+229) &\equiv0\pmod{7},\\
\Delta_{3}(343n+278) &\equiv0\pmod{7},\\
\Delta_{3}(343n+327) &\equiv0\pmod{7}.
\end{align*}

Quite recently, a variety of infinite families of congruences modulo 7 enjoyed by $\Delta_{3}(n)$ also have been found. See, for example, Xia \cite{Xia2015}, Yao and Wang \cite{YW2016}.

In this paper, we establish two infinite families of congruences modulo 5 and 25 for $\Delta_{k}(n)$ as follows:
\begin{theorem}\label{THM:mod 5}
For all $n\geq0$,
\begin{align}
\Delta_{k}(25n+24) &\equiv0\pmod{5}, \quad \emph{if}\quad k\equiv12\pmod{25},\label{con:mod 5}
\end{align}
\end{theorem}

\begin{theorem}\label{THM:mod 25}
For all $n\geq0$,
\begin{align}
\Delta_{k}(125n+99) &\equiv0\pmod{25}, \quad \emph{if}\quad k\equiv62\pmod{125}.\label{con:mod 25}
\end{align}
\end{theorem}

Moreover, we obtain the following infinite families of congruences modulo 7 and 49 for $\Delta_{k}(n)$.
\begin{theorem}\label{THM:mod 7}
For all $n\geq0$,
\begin{align}
\Delta_{k}(49n+s)\equiv0\pmod{7}, \quad \emph{if}\quad k\equiv24\pmod{49}\label{con:mod 7},
\end{align}
where $s=19, 33, 40$, and $47$.
\end{theorem}

\begin{theorem}\label{THM:mod 49}
\begin{align}
\Delta_{k}(343n+t)\equiv 0\pmod{49}, \quad \emph{if}\quad k\equiv171\pmod{343}.\label{con:mod 49}
\end{align}
where $t=96, 292$, and $341$.
\end{theorem}

\section{Proofs of Theorems \ref{THM:mod 5}--\ref{THM:mod 49}}
\subsection{The case mod 5}
To prove \eqref{con:mod 5}, we collect some useful identities. Recall that the Ramanujan theta function $f(a,b)$ is defined by
\begin{align}
f(a,b):=\sum_{n=-\infty}^{\infty}a^{n(n+1)/2}b^{n(n-1)/2},\label{Ramanujan theta func}
\end{align}
where $|ab|<1$. The Jacobi triple product identity can be stated as
\begin{align}
f(a,b)=(-a,-b,ab;ab)_{\infty}.\label{JTP}
\end{align}

One specialization of \eqref{Ramanujan theta func} is given by \cite{Ber}:
\begin{align*}
\psi(q) &:=f(q,q^3)=\sum_{n=0}^{\infty}q^{n(n+1)/2}.
\end{align*}
According to \eqref{JTP}, we have
\begin{align*}
\psi(q)=\dfrac{\f{2}^2}{\f{1}}.
\end{align*}

The following two identities are given by Hirschhorn \cite{Hir2017,Hir} and Berndt \cite{BerIII}.
\begin{lemma}
\begin{align}
\psi(q) &=a+qb+q^3c,\label{psi dissec 1}\\
\psi^{2}(q^{5}) &=ab+q^{5}c^{2}.\label{psi dissec 2}
\end{align}
where
\begin{align*}
a=f(q^{10},q^{15}), b=f(q^{5},q^{20}), c=\psi(q^{25}).
\end{align*}
\end{lemma}
Eq. \eqref{psi dissec 1} comes from \cite[Eq. (2.1)]{Hir2017}, \cite[Entry 10(i)]{BerIII} and Eq. \eqref{psi dissec 2} comes from \cite[Eq. (34.1.21)]{Hir}.

Employing the binomial theorem, we can easily establish the following congruence, which will be frequently used without explicit mention.
\begin{lemma}\label{cong lemma}
If $p$ is a prime, $\alpha$ is a positive integer, then
\begin{align*}
\f{\alpha}^{p} &\equiv\f{p\alpha}\pmod{p},\\
\f{1}^{p^{\alpha}} &\equiv\f{p}^{p^{\alpha-1}}\pmod{p^{\alpha}}.
\end{align*}
\end{lemma}

Now, we are ready to state the proof of \eqref{con:mod 5}.

Notice that $k\equiv12\pmod{25}$, then $2k+1\equiv0\pmod{25}$. Let $2k+1=25j$, where $j$ is a positive integer, we get, (all the following congruences are modulo 5)
\begin{align*}
\sum_{n=0}^{\infty}\Delta_{k}(n)q^{n} &=\dfrac{\f{2}\f{25j}}{\f{1}^{3}\f{50j}}=\dfrac{\psi^{3}(q)}{\f{2}^{5}}\dfrac{\f{25j}}{\f{50j}}\\
 &\equiv\dfrac{\psi^{3}(q)}{\f{10}}\dfrac{\f{25j}}{\f{50j}}.
\end{align*}

Invoking \eqref{psi dissec 1}, we obtain
\begin{align*}
\sum_{n=0}^{\infty}\Delta_{k}(n)q^{n} &\equiv\dfrac{\f{25j}}{\f{10}\f{50j}}\times\Bigg(a^{3}+3qa^{2}b+3q^{2}ab^{2}\\
 &\quad\quad+q^{3}b^{3}+3q^{3}a^{2}c+6q^{4}abc+3q^{5}b^{2}c+3q^{6}ac^{2}+3q^{7}bc^{2}+q^{9}c^{3}\Bigg).
\end{align*}

Extracting those terms involving the powers $q^{5n+4}$ and combining \eqref{psi dissec 2}, we have
\begin{align*}
\sum_{n=0}^{\infty}\Delta_{k}(5n+4)q^{5n} &\equiv\dfrac{\left(abc+q^{5}c^{3}\right)\f{25j}}{\f{10}\f{50j}}
 =\dfrac{c(ab+q^{5}c^{2})\f{25j}}{\f{10}\f{50j}}\\
 &=\psi(q^{25})\psi^{2}(q^{5})\dfrac{\f{25j}}{\f{10}\f{50j}},
\end{align*}
then
\begin{align*}
\sum_{n=0}^{\infty}\Delta_{k}(5n+4)q^{n} &\equiv\psi(q^{5})\psi^{2}(q)\dfrac{\f{5j}}{\f{2}\f{10j}}=\psi(q^{5})\dfrac{\f{2}^{3}}{\f{1}^{2}}
\dfrac{\f{5j}}{\f{10j}}\\
 &\equiv\psi(q^{5})\dfrac{\f{5j}}{\f{10j}}\dfrac{\f{1}^{3}\f{2}^{3}}{\f{5}}:=\sum_{n=0}^{\infty}a(n)q^{n},
\end{align*}
say. Thanks to the Jacobi's identity \cite[p. 14, Theorem 1.3.9]{Ber}
\begin{equation*}
\f{1}^3=\sum_{n=0}^\infty(-1)^n(2n+1)q^{n(n+1)/2},
\end{equation*}
it follows that
\begin{align}
\f{1}^3 &=J_{0}+J_{1}+J_{3},\nonumber\\
\f{2}^{3} &=J^{*}_{0}+J^{*}_{1}+J^{*}_{2},\label{2 JT}
\end{align}
where $J_{i}$ (resp. $J^{*}_{i}$) consists of those terms in which the power of $q$ is $i$ modulo 5. Furthermore, we see that $J_{3}\equiv 0\pmod{5}$ and $J^{*}_{1}\equiv0\pmod{5}$, so
\begin{align*}
(q;q)_\infty^3 &\equiv J_0+J_1 \pmod{5},\\
\f{2}^{3} &\equiv J^{*}_{0}+J^{*}_{2}\pmod{5}.
\end{align*}

Therefore,
\begin{align*}
\f{1}^3\f{2}^{3} \equiv\left(J_{0}+J_{1}\right)\left(J^{*}_{0}+J^{*}_{2}\right)\pmod{5},
\end{align*}
which contains no terms of the form $q^{5n+4}$. Hence $a(5n+4)\equiv0\pmod{5}$, equivalently,
\begin{align*}
\Delta_{k}(25n+24)\equiv0\pmod{5},
\end{align*}
where $k\equiv12\pmod{25}$, as desired.

To prove the remaining congruences \eqref{con:mod 25}--\eqref{con:mod 49}, we need to use a result of Radu and Sellers \cite[Lemma 2.4]{RS2011}. Before introducing the result of Radu and Sellers, we will briefly interpret some notations.

For a positive integer $M$, let $R(M)$ be the set of integer sequences $\{r:r=(r_{\delta_{1}},\cdots,r_{\delta_{k}})\}$ indexed by the positive divisors $1=\delta_{1}<\cdots<\delta_{k}=M$ of $M$.
For some $r\in R(M)$, define
$$f_r(q):=\prod_{\delta\mid M}(q^{\delta};q^{\delta})_{\infty}^{r_\delta}=\sum_{n=0}^{\infty}c_r(n)q^n.$$

Given a positive integer $m$. Let $\mathbb{Z}^{*}_{m}$ be the set of all invertible elements in $\mathbb{Z}_{m}$, and $\mathbb{S}_{m}$ be the set of all squares in $\mathbb{Z}^{*}_{m}$. We also define the set
$$P_{m,r}(t):=\left\{t'\ |\ t'\equiv ts+\frac{s-1}{24}\sum_{\delta\mid M}\delta r_{\delta}\pmod{m},0\le t'\le m-1,[s]_{24m}\in\mathbb{S}_{24m}\right\},$$
where $t \in \{0,\cdots,m-1\}$ and $[s]_{m}=s+m\mathbb{Z}$.

Let $\Gamma:=SL_2(\mathbb{Z})$ and $\Gamma_\infty:=\left\{\left.\begin{pmatrix}1 & h \\0 & 1 \end{pmatrix}\ \right|\ h\in\mathbb{Z}\right\}$. For a positive integer $N$, we define the congruence subgroup of level $N$ as
$$\Gamma_0(N):=\left\{\left.\begin{pmatrix}a & b\\c & d\end{pmatrix}\in\Gamma\ \right|\ c\equiv 0\pmod{N}\right\}.$$
The index of $\Gamma_0(N)$ in $\Gamma$ is given by
$$[\Gamma:\Gamma_0(N)]=N\prod_{p\mid N}(1+p^{-1}),$$
where the product runs through the distinct primes dividing $N$.

Let $\kappa=\kappa(m)=\gcd(m^2-1,24)$ and denote $\Delta^*$ by the set of tuples $(m,M,N,t,r=(r_\delta))$ satisfying conditions given in \cite[p. 2255]{RS2011}, we set
$$p_{m,r}(\gamma)=\min_{\lambda\in\{0,\ldots,m-1\}}\frac{1}{24}\sum_{\delta\mid M}r_\delta\frac{\gcd^2(\delta(a+\kappa\lambda c),mc)}{\delta m}$$
and
$$p_{r'}^*(\gamma)=\frac{1}{24}\sum_{\delta\mid N} \frac{r'_\delta\gcd^2(\delta,c)}{\delta},$$
where $\gamma=\begin{pmatrix}a & b \\c & d \end{pmatrix}$, $r\in R(M)$, and $r'\in R(N)$.

The lemma of Radu and Sellers is stated as follows.
\begin{lemma}\label{le:RS}
Let $u$ be a positive integer, $(m,M,N,t,r=(r_\delta))\in\Delta^*$, $r'=(r'_\delta)\in R(N)$, $n$ be the number of double cosets in $\Gamma_0(N)\backslash\Gamma/\Gamma_\infty$ and $\{\gamma_1,\ldots,\gamma_n\}$ $\subset\Gamma$ be a complete set of representatives of the double coset $\Gamma_0(N)\backslash\Gamma/\Gamma_\infty$. Assume that $p_{m,r}(\gamma_i)+p_{r'}^*(\gamma_i)\ge 0$ for all $i=1,\ldots,n$. Let $t_{\min} := \min_{t'\in P_{m,r}(t)}t'$ and
$$v:=\frac{1}{24}\left(\left(\sum_{\delta\mid M}r_\delta+\sum_{\delta\mid N}r'_\delta\right)[\Gamma:\Gamma_0(N)]-\sum_{\delta\mid N}\delta r'_\delta\right)-\frac{1}{24m}\sum_{\delta\mid M}\delta r_\delta-\frac{t_{\min}}{m}.$$
Then if
$$\sum_{n=0}^{\lfloor v \rfloor}c_r(mn+t')q^n\equiv 0 \pmod{u}$$
for all $t'\in P_{m,r}(t)$, then
$$\sum_{n=0}^{\infty}c_r(mn+t')q^n\equiv 0 \pmod{u}$$
for all $t'\in P_{m,r}(t)$.
\end{lemma}

\subsection{The case mod 25}\label{subsec:mod 25}
Let
\begin{align*}
\sum_{n=0}^{\infty}b(n)q^n=\frac{\f{2}}{\f{1}^3}.
\end{align*}
By Lemma \ref{cong lemma}, we obtain
\begin{align}
\sum_{n=0}^{\infty}b(n)q^n =\frac{\f{2}}{\f{1}^3}\equiv \frac{\f{1}^{22}\f{2}}{\f{5}^5}\pmod{25}.\label{cong: mod 25 equlity}
\end{align}
In this case, we may take
$$(m,M,N,t,r=(r_1,r_2,r_{5},r_{10}))=(125,10,10,99,(22,1,-5,0))\in\Delta^*.$$
By the definition of $P_{m,r}(t)$, we have $P_{m,r}(t)=\{99\}$.
Now we can choose
$$r'=(r'_1,r'_2,r'_5,r'_{10})=(13, 0, 0, 0).$$

Let
$$\gamma_\delta=\begin{pmatrix}
1 & 0 \\
\delta & 1
\end{pmatrix}.$$
Radu and Sellers \cite[Lemma 2.6]{RS2011} also proved that $\{\gamma_\delta:\delta\mid N\}$ contains a complete set of representatives of the double coset $\Gamma_0(N)\backslash\Gamma/\Gamma_\infty$.

One readily verifies that all assumptions of Lemma \ref{le:RS} are satisfied. Furthermore we obtain the upper bound $\lfloor v \rfloor=21$.

For all congruences in \eqref{cong:mod 25}, we check that they hold for $n$ from 0 to their corresponding upper bound $\lfloor v \rfloor$ via \emph{Mathematica}. It follows by Lemma \ref{le:RS} and \eqref{cong: mod 25 equlity} that
\begin{align}
b(125n+99) &\equiv0\pmod{25}\label{cong:mod 25}
\end{align}
holds for all $n\geq0$. When $k\equiv62\pmod{125}$, i.e., $2k+1\equiv0\pmod{125}$, \eqref{con:mod 25} is a direct consequence of \eqref{cong:mod 25}.

\subsection{The case mod 7}
Firstly, we have
\begin{align*}
\sum_{n=0}^{\infty}b(n)q^n=\frac{\f{2}}{\f{1}^3}\equiv \frac{\f{1}^{4}\f{2}}{\f{7}}\pmod{7}.
\end{align*}
To prove \eqref{con:mod 7}, it suffices to show
\begin{align}
b(49n+s) &\equiv0\pmod{7}\label{cong: mod 7}
\end{align}
for $s=19, 33, 40$, and 47.

We first show the cases $s=19, 33$, and 40. Taking
$$(m,M,N,t,r=(r_1,r_2,r_{7},r_{14}))=(49,14,14,33,(4,1,-1,0))\in\Delta^*.$$
We compute that $P_{m,r}(t)=\{19,33,40\}$. Now we can choose
$$r'=(r'_1,r'_2,r'_7,r'_{14})=(3, 0, 0, 0).$$
and taking $\gamma$ as in Subsection \ref{subsec:mod 25}, we verify that all these constants satisfy the assumption of Lemma \ref{le:RS}. We thus obtain $\lfloor v \rfloor=6$. With the help of \emph{Mathematica}, we see that \eqref{cong: mod 7} holds up to the bound $\lfloor v \rfloor$ with $t\in\{19, 33, 40\}$, and therefore it holds for all $n\geq0$ by Lemma \ref{le:RS}.

Now we will turn to the case $s=47$. Again we take
\begin{align*}
(m,M,N,t,r=(r_1,r_2,r_{7},r_{14}))=(49,14,14,47,(4,1,-1,0))\in\Delta^*
\end{align*}
and
\begin{align*}
r'=(r'_1,r'_2,r'_7,r'_{14})=(3, 0, 0, 0).
\end{align*}
In this case we obtain $P_{m,r}(t)=\{47\}$. One readily computes that $\lfloor v \rfloor=6$. Thus we verify the first 6 terms of \eqref{cong: mod 7} via \emph{Mathematica}. It follows by Lemma \ref{le:RS} that it holds for all $n\geq0$.

\subsection{The case mod 49}
Similarly, we have
\begin{align*}
\sum_{n=0}^{\infty}b(n)q^n=\frac{\f{2}}{\f{1}^3}\equiv \frac{\f{1}^{46}\f{2}}{\f{7}^{7}}\pmod{49}.
\end{align*}
To prove \eqref{con:mod 49}, it needs to show
\begin{align}
b(343n+t) &\equiv0\pmod{49}\label{cong: mod 49}
\end{align}
for $t=96$, 292, and 341.

In these cases, we may set
$$(m,M,N,t,r=(r_1,r_2,r_{7},r_{14}))=(343,14,14,96,(46,1,-7,0))\in\Delta^*.$$
We compute that $P_{m,r}(t)=\{96, 292, 341\}$. Now we can choose
$$r'=(r'_1,r'_2,r'_7,r'_{14})=(18, 0, 0, 0).$$
We thus obtain $\lfloor v \rfloor=56$. Similarly we verify the first 56 terms of \eqref{cong: mod 49} through \emph{Mathematica}. It follows by Lemma \ref{le:RS} that
\begin{align*}
\Delta_{k}(343n+96)\equiv\Delta_{k}(343n+292)\equiv\Delta_{k}(343n+341)\equiv 0\pmod{49}
\end{align*}
holds for all $n\geq0$ if $k\equiv171\pmod{343}$.

\section{Final remarks}
Despite the universality of Radu and Sellers' lemma, we should point out that their method is not elementary because it highly relies on modular forms. On the other hand, the proofs of \eqref{con:mod 25}--\eqref{con:mod 49} are routine to some extent. It would be interesting to find elementary proofs of these congruences.

\section*{Acknowledgement}
The author would like to thank Shishuo Fu, Shane Chern and Michael D. Hirschhorn for their helpful comments and suggestions on the original manuscript. This work was supported by the National Natural Science Foundation of China (No.~11501061).


\begin{thebibliography}{99}


\bibitem{AP2007}G.~E. Andrews and P. Paule, MacMahon's partition analysis XI, broken diamonds and modular forms, \textit{Acta Arith.} {\bf 126} (3) (2007) 281--294.

\bibitem{Ber} B. C. Berndt, \textit{Number theory in the spirit of Ramanujan}, Student Mathematical Library, \textbf{34}. American Mathematical Society, Providence, RI, 2006. xx+187 pp.

\bibitem{BerIII} B. C. Berndt, \textit{Ramanujan Notebooks Part III}, Springer--Verlag, New York Inc., 1991. xx+510 pp. 

\bibitem{Chan2008}S. H. Chan, Some congruences for Andrews-Paule's broken 2-diamond partitions, \textit{Discrete Math.} \textbf{308} (2008) 5735--5741.

\bibitem{Hir2017}M. D. Hirschhorn, Broken 2-diamond partitions modulo 5, \textit{Ramanujan J.} to appear.

\bibitem{Hir}M. D. Hirschhorn, \textit{The power of $q$}, Developments in Mathematics Vol. 49, Springer 2017. xx+413 pp.

\bibitem{Jam2013}M. Jameson, Congruences for broken $k$-diamond partitions, \textit{Ann. Comb.} \textbf{17} (2013) 333--338.

\bibitem{PR2010}P. Paule and C. S. Radu, Infinite families of strange partition congruences for broken 2-diamonds, \textit{Ramanujan. J} \textbf{23} (2010) 409--416.

\bibitem{Rad2015}C. S. Radu, An algorithmic approach to Ramanujan-Kolberg identities, \textit{J. Symb. Comput.} \textbf{68} (1) (2015) 225--253.

\bibitem{RS2011}C. S. Radu and J. A. Sellers, Congruence properties modulo $5$ and $7$ for the $\mathrm{pod}$ function, \textit{Int. J. Number Theory} \textbf{7} (2011), no. 8, 2249--2259.

\bibitem{Xia2015}Ernest X. W. Xia, Infinite families of congruences modulo 7 for broken 3-diamond partitions, \textit{Ramanujan J.} to appear.

\bibitem{Xia2017}Ernest X. W. Xia, More infinite families of congruences modulo 5 for broken 2-diamond partitions, \textit{J. Number Theory.} \textbf{170} (2017) 250--262.

\bibitem{YW2016}Olivia X. M. Yao and Y. J. Wang, Newman's identity and infinite families of congruences modulo 7 for broken 3-diamond partitions, \textit{Ramanujan J.} to appear.

\end{thebibliography}
\end{document}